\documentclass[10pt,a4paper]{article}
\usepackage[utf8]{inputenc}
\usepackage{amsmath}
\usepackage{amsfonts}
\usepackage{amsthm}
\usepackage{calrsfs}
\usepackage{amssymb}
\usepackage{xcolor}
\usepackage{dsfont}
\usepackage{mathtools}

\usepackage{geometry}
\theoremstyle{plain}
\newtheorem{theorem}{Théorème}
\newtheorem{prop}[theorem]{\textbf{Proposition}}
\newenvironment{dem}{\textbf{Proof.}}{\qed \\ }
\newenvironment{demof}{\textbf{Proof of }}{\qed \\ }
\newtheorem{corol}[theorem]{\textbf{Corollary}}
\newtheorem{lemm}[theorem]{\textbf{Lemma}}
\newtheorem{thm}[theorem]{\textbf{Theorem}}

\newtheorem*{assumpt}{\textbf{Assumptions on $(L,\mu)$}}
\newtheorem*{assum}{\textbf{Assumption 3}}

\newtheorem{defi}{\textbf{Definition}}

\usepackage{filecontents,lipsum}

\usepackage{mystyle}

\title{Stability of the Poincaré constant}
\author{Jordan Serres}

\begin{document}

\maketitle

\begin{abstract} 
We study stability of the sharp Poincaré constant of the invariant probability measure of a reversible diffusion process satisfying some natural conditions. The proof is based on the spectral interpretation of Poincaré inequalities and Stein's method. In particular, these results are applied to the gamma distributions, to the Brownian motion on spheres and to the Brascamp-Lieb inequality for one-dimensional log-concave measures.
\end{abstract}

\section{Introduction}

A probability measure $\mu$ on $\mathbb{R}^d$ is said to satisfy a Poincaré inequality when there exists a positive finite constant $C$ such that for all functions $f$ in the Sobolev space $H^{1}(\mu)$,
\begin{equation}
\mathrm{Var}_\mu (f) \leq C \int |\nabla f|^2 d\mu.  
\end{equation}
We denote by $C_P(\mu)$ the smallest constant for which the above inequality holds.
Poincaré inequalities have many applications (see for instance the survey \cite{surveyineg}). For instance they can be seen as embeddings of weighted Sobolev spaces in $L^2$, or as quantifying concentration of measure phenomenon (see e.g. \cite{Troyanov}). The sharp Poincaré constant governs ergodicity of the underlying dynamic in the $L^2$ sense, and the convergence rate of some algorithms used for numerical simulations (see e.g. \cite{kls, klssurvey}). When $\mu$ is reversible for a Markov process, a Poincaré inequality has a spectral interpretation: the infinitesimal generator $L$ of the Markov process is symmetric on $L^2(\mu)$ and the quantity $\lambda_\mu := \frac{1}{C_P(\mu)}$ is then the spectral gap of the positive symmetric operator $-L$ (see \cite[section 4.2.1]{BGL}).

Stability results for Poincaré constant began to appear in the late 80's. Chen \cite[Corollary 2.1]{1987} showed that all isotropic probability measures on $\mathbb{R}^d$ have sharp Poincaré constant greater than $1$. He proved furthermore that the standard Gaussian is the only one attaining $1$. Then Utev \cite{U} refined this result in dimension one, quantifying the difference between Poincaré constants in term of total variation distance: $$C_P(\nu)\geq 1 + \frac{1}{9}d_{TV}(\nu,\gamma)^2 $$ where $\nu$ is a normalized probability measure on $\mathbb{R}$, $\gamma$ is the standard Gaussian and $d_{TV}$ is the total variation distance.

More recently, Courtade, Fathi and Pananjady \cite{existsteinkernel}, extended it to the multidimensional case with the Wasserstein-$2$ distance:
\begin{equation}\label{stabconnu}
C_P(\nu)\geq 1 + \frac{W_2(\nu,\gamma)^2}{d} 
\end{equation}
where $\nu$ is a centered probability measure on $\mathbb{R}^d$, normalized such that $\int|x|^2\,d\nu=d$, $\gamma$ denotes the Gaussian $\mathcal{N}(0,I_d)$ and $W_2$ is the $2$-Wasserstein distance (see \cite[chapter 6]{villani}). Our goal is to get such stability results in a more abstract framework, say for a general reference probability measure $\mu$ on a manifold instead of the Gaussian on $\mathbb{R}^d$. What we call stability results for Poincaré constant with respect to some distance $d$, are inequations of the form
\begin{equation}\label{stabdef}
 d(\mu,\nu) \leq \phi(C_P(\mu),C_P(\nu))
\end{equation} where $d$ is a distance on the space of probability measures, $\phi: \{(x,y)\in \mathbb{R}^2\,|\,y\geq x>0\}\rightarrow\mathbb{R}_+$ is a continuous function such that $\forall x>0$, $\phi(x,x)=0$, and $C_P(\mu)$, $C_P(\nu)$ are respectively sharp Poincaré constants for $\mu$ a reference measure, and $\nu$ satisfying some constraints.

From now on, let us consider $L$ an infinitesimal reversible Markov diffusion generator on a Riemannian manifold $M$ with domain $\mathcal{D}$, see \cite[Sect. 1.4]{BGL} for further details. Let $\Gamma(f,g):=\frac{1}{2}\left(L(fg)-fLg-gLf\right)$ be the carré du champ operator. To lighten the notation, we set $\Gamma(f):=\Gamma(f,f) $. The diffusion property means that $\forall \phi\in C^2(\mathbb{R})$, $\Gamma(\phi(f))=\phi'(f)^2 \Gamma (f)$ and also $ L(\phi(f))=\phi''(f)\Gamma (f) + \phi'(f)L(f) $, see \cite[Section 1.11]{BGL}. Let us assume $\mu$ to be the only reversible probability measure on $M$ for the Markovian process generated by $L$. This existence and uniqueness assumption is verified for instance if the process is irreducible and strongly Feller (see \cite[Chapter 4]{Khasminskii}).
 
Assume that $\mu$ satisfies the following Poincaré inequality:
\begin{equation}
\forall f\in H^1(\mu), \, \mathrm{Var}_\mu (f) \leq C_P(\mu) \int \Gamma(f)\, d\mu  
\end{equation}
 where 
 \begin{equation}\label{defH1}
 H^1(\mu):=\{f\in L^2(\mu)\, |\ \int f\, d\mu = 0 \,\, \mathrm{and}\, \int\Gamma (f) d\mu<\infty\} 
 \end{equation} and $C_P(\mu)=\underset{f\in H^1(\mu)}{\sup} \frac{\mathrm{Var}_\mu (f)}{\int\Gamma (f) d\mu}<\infty$ is the sharp Poincaré constant.

\begin{assumpt} \,\newline

\begin{enumerate}

\item \label{spectralgapattained} There exists an eigenfunction $f_0$ attaining the spectral gap: $-Lf_0=\lambda_\mu f_0$. We will always choose it to be normalized with respect to $\mu$, i.e. such that $\int f_0d\mu=0$ and $\int f_0^2d\mu=1 $.

\item \label{gammaf0measuable} $\Gamma(f_0)$ can be written as $h\circ f_0$ for some smooth function $h:I\rightarrow \mathbb{R}_+$ with $I:=f_0(M)$.

\item \label{criticalmaxima} 
The function $h$ does not vanish on the interior of $I$.

\end{enumerate}
\end{assumpt} 
Assumption $\ref{spectralgapattained}$ is used to push forward the diffusion from $M$ to $I\subset \mathbb{R}$ through $f_0$ (see Section $\ref{quotienprocess}$). This assumption is necessary. Indeed for $M=\mathbb{R}$ and $Lf(x)=f''(x)-\mathrm{sgn}(x)f'(x)$ where $\mathrm{sgn}(x)$ denotes the sign of $x$, the spectral gap is not attained (see \cite[Section 4.4.1]{BGL}). %
Assumption $\ref{gammaf0measuable}$ guarantees the push forward process to be Markovian (see Section $\ref{quotienprocess}$). This assumption would obviously be verified as soon as $f_0$ is injective. In dimension $1$, it is always the case (see Lemma $\ref{monotoniefirstfp}$).
Assumption $\ref{criticalmaxima}$ guarantees the diffusion on $I$ to be irreducible (see Section $\ref{nonvanish}$).\\

Let $\nu$ be another measure on $M$ satisfying the same Poincaré inequality with sharp constant $C_P(\nu)$: $$ \forall f\in H^1(\nu), \ \mathrm{Var}_\nu (f) \leq C_P(\nu) \int\Gamma (f) d\nu.  $$ where $H^1(\nu)$ is defined similarly as above. We also ask $\nu$ to be $f_0$-normalized, that is $f_0\in H^1(\nu)$ and
\begin{equation}\label{hypnu}
\int f_0\,d\nu=0, \ \int f_0^2 d\nu=1,\ \mathrm{ and} \int\Gamma (f_0)\,d\nu \leq \frac{1}{C_P(\mu)}. 
\end{equation}
The class of measures satisfying the three equations $\int f_0\,d\nu=0, \ \int f_0^2 d\nu=1,\ \mathrm{ and} \int\Gamma (f_0)\,d\nu = \frac{1}{C_P(\mu)}$ being a space of codimension $3$ in the space of all probability measures on $M$, one sees that ($\ref{hypnu}$) is satisfied on a half-space of codimension $3$ in the infinite dimension space of probability measures.

Applying Poincaré inequality for $\nu$ to $f_0$, one immediately gets $C_P(\nu)\geq C_P(\mu)$.
We refine this minorization by proving the following stability theorem:
\begin{thm}\label{thmprincipal}
Under Assumptions $\ref{spectralgapattained}$, $\ref{gammaf0measuable}$ and $\ref{criticalmaxima}$,
\begin{equation}\label{stabeq}
W_1(\mu^*,\nu^*) \leq C_h\left(\sqrt{\delta} + \sqrt{C_P(\nu)}\,\delta \right),
\end{equation}
where $\mu^*$ (resp. $\nu^*$) is the pushforward of $\mu$ (resp. $\nu$) by $f_0$, $W_1$ is the $1$-Wasserstein distance (see Definition $\ref{w1int}$), $\delta:= \frac{ C_P(\nu) - C_P(\mu)}{C_P(\nu)\,C_P(\mu)}$, and $C_h$ is a constant defined in Proposition $\ref{gammasteinbound}$ whose finitness only depends on the behavior of $h$ at the boundary of $I$.
\end{thm} %

The method followed in this article is to derive some approximate integration by part formula for $\nu$ (see Section $\ref{sectipp}$) and then to use Stein's method (see Section $\ref{steinmethod}$). The problem is then reduced to bound the carré du champ operator of the solution (see Section $\ref{steinsolsection}$). \\

A particuliar case, treated by Courtade and Fathi \cite{existsteinkernel,CF}, is when $L=\Delta -x\cdot \nabla$ is the Ornstein-Uhlenbeck process on $\mathbb{R}^d$. The invariant measure is the standard Gaussian measure $\gamma$ and the $d$ coordinate projections $x_1,...x_d$ are orthogonal eigenfunctions satisfying Assumption $\ref{spectralgapattained}$, $\ref{gammaf0measuable}$ and $\ref{criticalmaxima}$. Choosing any of these projections gives the one dimensional Gaussian $\mathcal{N}(0,1)$ for the pushforward measure. In this case, the approach described above works using only classical results of Stein's method for the Gaussian (see \cite[Section 2]{fundamStein}). Moreover, if we use the same approach with $f_0=(x_1,...,x_d)$ a vector valued function, it works thanks to technical bounds obtained in \cite[Lemma 3.3]{multivariate} and it gives similar results as $\eqref{stabconnu}$. Taking vector valued eigenfunction when the dimension of the eigenspace is greater than $1$ appears to be the natural extension of this method. Technical issues come here from the fact that Stein's method involves control of the Hessian of solutions to a family of PDE.\\

In Section $\ref{gammadistrib}$, we obtain stability result for the Laguerre process as Corollary of Theorem $\ref{thmprincipal}$. As far as we know, stability results had not yet been investigated in this case.\\

One can ask for another distance than the $1$-Wasserstein in the stability inequality ($\ref{stabeq}$). The main difference lies in the set of target function chosen using Stein's method. We will see in Section $\ref{ellipt}$ that replacing the condition $C_h<\infty$ by $h>\kappa>0$, where $\kappa$ is a constant, gives a stability result in total variation distance (and hence in Kolmogorov distance). We illustrate this result in Section $\ref{brascampliebstab}$ with uniformly log-concave measures.\\

A natural question arising is then to compare stability results involving different distances. We will say that a $d_1$ stability result 
\begin{equation}\label{d_1}
 d_1(\mu,\nu) \leq \phi_1(C_P(\mu),C_P(\nu)) 
\end{equation} is stronger than a $d_2$ one 
\begin{equation}\label{d_2}
 d_2(\mu,\nu) \leq \phi_2(C_P(\mu),C_P(\nu)) 
\end{equation} if ($\ref{d_1}$) implies the existence of $\phi_2$ in ($\ref{d_2}$). Comparison between stability results in Kolmogorov distance obtained from Theorems $\ref{stabilityhminore}$ and $\ref{thmcentralexploitable}$ shall be discussed in Section $\ref{comparison}$.

\section{The quotient process}\label{quotienprocess}

In the one dimensional Ornstein-Uhlenbeck process, the first non zero eigenfunction is injective. So one may expect the important properties of $L$ to be preserved when mapping it onto $I$ through such eigenfunction.
In a broader context, according to \cite[page 60]{BGL}, Assumption $\ref{gammaf0measuable}$ implies that all of the Markovian structure on the manifold is mapped onto a Markovian structure on $I:=f_0(M)\subset\mathbb{R}$. Assuming that $M$ is connected, $I$ is an interval since $f_0$ is continuous. We define $a:=\inf I$ and $b:= \sup I$. Moreover, the diffusion property allows to write that $$\forall \phi\in C^2(\mathbb{R}),\quad L(\phi(f_0)) = \phi''(f_0)\Gamma f_0+ \phi'(f_0)Lf_0  =  \phi''(f_0)h(f_0)  - \phi'(f_0)\frac{1}{C_{P}(\mu)}f_0.$$
Hence, the induced Markov process has generator
\begin{equation}\label{lstardirect}
L^*(\phi)(x) := h(x)\phi''(x) - \frac{x}{C_{P}(\mu)}\phi'(x) 
\end{equation}
and reversible measure $\mu^*:=f_0^{\#}\mu$. This one dimensional Markov process will be called the quotient process.

If $h$ is constant, then it is equal to $\frac{1}{C_P(\mu)}$. Indeed, one gets that $\Gamma(f_0) = h$, but we know that $\int \Gamma(f_0) d\mu = \frac{1}{C_P(\mu)}$. In this case, $L^*$ is reduced to the Ornstein-Uhlenbeck generator with reversible measure $\mathcal{N}(0,\frac{1}{C_P(\mu)})$. Therefore, $\mu^* = \mathcal{N}(0,\frac{1}{C_P(\mu)})$.  Similarly, $\nu^*:=f_0^{\#}\nu$ will satisfy Poincaré inequality with sharp constant $C_{P}(\nu^*)\leq \frac{C_P(\nu)}{C_P(\mu)}$. Indeed, using Poincaré inequality for $\nu$, $$\mathrm{var}_{\nu^*}(\phi) = \mathrm{var}_{\nu}(\phi\circ f_0)\leq C_{P}(\nu)\int \Gamma(\phi\circ f_0)d\nu =  C_{P}(\nu)\int \phi'(f_0)^2 h(f_0)d\nu = \frac{C_P(\nu)}{C_P(\mu)}\int \phi'^2\,d\nu^*.$$ Of course the carré du champ operator used in Poincaré inequalities on $\mathbb{R}$ is the square of the first derivative. 
At that point, we can apply known stability results (see \cite{existsteinkernel}) and obtain: $$ \frac{C_P(\nu)}{C_P(\mu)}\geq C_{P}(\nu^*)\geq C_{P}(\mu) + W_{2}(\nu^*, \mu^*)^2. $$ 
In the sequel, the goal will be to prove such stability inequalities in a broader context, that is without assuming $h$ to be constant.
Let us now compute the carré du champ of the quotient process.
\begin{prop}\label{poincarpourmustar}
The carré du champ operator associated to $(L^*,\mu^*)$ is $$ (\Gamma^*\psi) (t) = h(t)\psi'(t)^2.$$
Moreover, with this operator, $\mu^*$ satisfies a Poincaré inequality with sharp constant $C_P(\mu)$ and the inequality becomes an equality for $\psi=Id$.
\end{prop}
\begin{dem}
It is simply a computation using the basic definition of $\Gamma f := \frac{1}{2}[L(f^2)-2fLf]$. The Poincaré inequality for $\mu^*$ follows from the Poincaré inequality for $\mu$: $$\mathrm{var}_{\mu^*}(\psi) = \mathrm{var}_{\mu}(\psi\circ f_0)\leq C_{P}(\mu)\int \Gamma(\psi\circ f_0)d\mu =  C_P(\mu)\int \psi'(f_0)^2 h(f_0)d\mu =C_P(\mu)\int \Gamma^*\psi d\mu^*.$$
 Now taking $\psi=Id$, the inequality becomes an equality because of the definition of $f_0$, showing that $C_P(\mu)$ is sharp.
\end{dem}
Replacing $\mu^*$ by $\nu^*$ in the proof above, we get that $\nu^*$ statisfies a Poincaré inequality with constant $C_P(\nu)$. However in this case, $C_P(\nu)$ may not be the sharp constant for $\nu^*$.

\subsection{The non vanishing assumption}\label{nonvanish}

Let us consider the case where $Lf:= \alpha f''+\beta f' $ is a diffusion operator on an interval $M=J:=(c,d)\subset	\mathbb{R}$, with $\alpha$ and $\beta$ two continuous functions on $J$ such that $\alpha > 0$ and $$\rho(x):=\frac{1}{\alpha(x)}\exp\left(\int_c^x \frac{\beta(t)}{\alpha(t)}\,dt \right)$$ is well defined on $J$. So $d\mu:= \rho(x)dx$ is reversible for $L$. Assume the operator $L$ satisfies both a Poincaré inequality with sharp constant $C_P$ and Assumption $\ref{spectralgapattained}$. Then the following holds.
\begin{lemm}\label{monotoniefirstfp}
The first eigenfunction $f_0$ is strictly monotone.
\end{lemm}
\begin{dem}
We extend the proof in e.g. \cite[Section 2]{poincaroninterval} which treats the case $\alpha\equiv 1$.
First, recall that $f_0$ is a minimizer of the Rayleigh ratio $\frac{\int_J \Gamma(f_0)d\mu}{\mathrm{var}_\mu (f_0)}$, and introduce $g(x) :=  \int_c^x |f'_0(t)|\,dt $. Then $g'(x)= |f'_0(x)| $ hence $\Gamma(g)=\Gamma(f_0)$. So $g$ has same energy than the eigenfunction $f_0$: $$\int_J \Gamma(g)\,d\mu = \int_J \Gamma(f_0)\,d\mu $$
On the other hand, $g$ has greater variance than $f_0$. Indeed, $$\mathrm{Var}_\mu(g) = \frac{1}{2}\int_J \int_J \left(\int_x^y |f_0'(t)|\,dt \right)^2d\mu(x)d\mu(y) \geq \frac{1}{2}\int_J \int_J \left(\int_x^y f_0'(t)\,dt \right)^2d\mu(x)d\mu(y) = \mathrm{Var}_\mu(f_0)$$ But $f_0$ being a minimizer of the Rayleigh ratio, one can infer, since $f_0 \in C^1(J)$, that $f'_0$ has same sign on the interval $J$, hence $f_0$ is monotone.
Assume now $f_0'\geq 0$ and let us show that $f_0'>0$. The eigenfunction $f_0$ being centered and continuous, it is negative in a neighborhood of $c$, then it is zero, and then positive up to $d$.
Using classical ODE tools for the equation $\alpha f''_0 + \beta f'_0 = \frac{-1}{C_P}f_0 $, one gets the formulas:
\begin{align}
 f_0'(x) &=\label{droite} \frac{1}{C_P\,\alpha(x)\rho(x)}\int_x^d f_0(t)d\mu(t) \\
 &=\label{gauche} \frac{-1}{C_P\,\alpha(x)\rho(x)}\int_c^x f_0(t)d\mu(t) 
\end{align}
Hence if $f_0'(x)=0$ with $x\in (c,d)$, then either $f_0(x)<0$ but Formula $\ref{gauche}$ implies $f_0'(x)>0$, or $f_0(x)>0$ but Formula $\ref{droite}$ implies $f_0'(x)>0$, or $f_0(x)=0$ but both formulas give then again $f_0'(x)>0$. 
 
Hence $f'_0$ is positive on $J$, justifying the claim that $f_0$ is strictly monotone.
\end{dem}

In this case the monotonicity implies that $h$ cannot vanish in the interior of $I$. This motivates Assumption $\ref{criticalmaxima}$ already mentioned in the introduction. Furthermore, since $a$ (respectively $b$) is the global minimum (resp. maximum) of $f_0$, it can be reformulated as:
\begin{assum}
All $x_0\in M$ such that $\Gamma(f_0)(x_0) = 0$ are global extrema of $f_0$.
\end{assum}

If $\Gamma=|\nabla|^2$ for the metric of $M$, then the assumption can be reformulated: all critical points of $f_0$ are global extrema. The following is immediate and justify why we often write $I$ instead of its interior $\overset{\circ}{I}$ in the sequel.

\begin{prop}
The eigenfunction $f_0$ satisfies Assumption $\ref{criticalmaxima}$ if, and only if, $h$ does not vanish on $\overset{\circ}{I}$.
\end{prop}

\begin{dem}
Assume that Assumption $\ref{criticalmaxima}$ holds, and let $t\in I$ such that $h(t)=0$. Then there exists $x\in M$ such that $t=f_0(x)$ and so $\Gamma(f_0(x))=0 $. Hence $x$ is a global extremum of $f_0$. Hence $t=a\notin\overset{\circ}{I}$ or $t=b\notin\overset{\circ}{I}$.
Conversely let $x\in M$ such that $\Gamma (f_0)(x)=0$. Then $f_0(x)$ is a zero of $h$. Hence by assumption it is on the boundary of $I$, so it is a global extremum of $f_0$.
\end{dem}

\subsection{Density of $\mu^*$}\label{zeroestdansI}

First, we point out that $\int f_0\,d\mu \in I$. Indeed $f_0$ is continuous, so $I=f_0(M)$ is connected in $\mathbb{R}$, hence convex. Therefore its expectation belongs to $I$. Recall that $\lambda_\mu:=\frac{1}{C_P(\mu)}$.
We now define the function $v : I\rightarrow \mathbb{R}_+$ by 
\begin{equation}\label{defdev}
v(t) := \exp\left( -\lambda_\mu\int_0^t \frac{u\,du}{h(u)} \right) 
\end{equation} 
Taking $0$ as reference point in the integral is justified by the above observation since $f_0$ was choosen to be centered. Moreover Assumption $\ref{criticalmaxima}$ garuantees that dividing by $h$ makes sense. 
\begin{prop}\label{densitymustar}
The measure $$ \rho(x)dx:=\frac{1}{ h(x)}\exp\left( -\lambda_\mu\int_{0}^{x}\frac{u\,du}{ h(u)}\right) \, \textbf{1}_{I} dx $$ is invariant for the quotient process, where $\textbf{1}$ denotes an indicator.
\end{prop}

\begin{dem}
It is enough to show that $$  \int_I L^* (f) \frac{1}{ h(x)}\exp\left( -\lambda_\mu \int_{0}^{x}\frac{u\,du}{h(u)}\right) dx = 0$$ for all functions $f$ compactly supported in $I$. Let $f$ be such a function. Recall that $L^* f = h f'' -\frac{x}{C_P(\mu)} f'$. We first compute by an integration by parts: $$ \int_I f''(x)v(x)dx = \underset{=0}{\underbrace{\left[f'(x)v(x)\right]_{a}^b}} + \lambda_\mu \int_I \frac{xf'(x)}{h(x)}v(x)dx,$$ and this allows us to conclude. 
\end{dem}
We know that the quotient process admits $f_0^{\#}\mu$ as invariant probability measure. Since we assumed $\mu$ to be the only invariant probability measure for $L$, Proposition $\ref{densitymustar}$ gives that the measure $\rho(x)dx$ is finite and then $f_0^{\#}\mu=\frac{1}{Z}\rho(x)dx$ where $Z$ is a normalization constant.
Let us conclude this section with a control of the tail of $\mu^*$. We denote the cumulative distribution of $\mu^*$ by $$q(t):=\mu^*(]-\infty,t])=\int_{-\infty}^{t}\frac{\mathbf{1}_I}{Z h(y)}\exp\left( - \lambda_\mu\int_{0}^{y}\frac{u\,du}{ h(u)}\right)  dy.$$ Recall that $a:=\inf I$, $b:= \sup I$ and let us denote $I_-:=I\cap\mathbb{R}_-$ and $I_+:=I\cap\mathbb{R}_+$.
\begin{lemm}\label{tailcontrol}
One can bound the tail $q$ of $\mu^*$ as follows:
\begin{itemize}
\item $\forall t\in I_-$, $$q(t)\leq \min\left(q(0),\frac{-C_P(\mu)}{Zt} \right)\,v(t).$$
\item $\forall t\in I_+$, $$1-q(t)\leq \min\left((1-q(0)), \frac{C_P(\mu)}{Zt}\right)\, v(t).  $$
\end{itemize}
\end{lemm}

\begin{dem}
We only process on $I_{-}$. The proof is similar on $I_{+}$.
\begin{enumerate}
\item Let us define $f : I_-\rightarrow \mathbb{R}$, by $f(t)=q(t)-q(0)v(t)$. We want to show that $f\leq 0$. Compute $$f'(t)= \frac{1}{Z h(t)}\exp\left( -\lambda_\mu \int_{a}^{t}\frac{u\,du}{ h(u)}\right) +\lambda_\mu q(0)\frac{t}{ h(t)}v(t) =\frac{v(t)}{h(t)}\left( \frac{1}{Z} +\lambda_\mu q(0)t \right). $$ If $-\frac{C_P(\mu)}{q(0)\,Z}<a$, then $f'\geq 0$ so $f\leq f(0)=0 $.\\
Else, $f$ decreases on $(a,-\frac{C_P(\mu)}{q(0)\,Z})$ and increases on $(-\frac{C_P(\mu)}{q(0)\,Z},0)$ but $\underset{t\rightarrow a}{\lim}f(t) \leq 0$ and $f(0)=0$, so we have the first claim.

\item  Let us define $f : I_-\rightarrow \mathbb{R}$, by $f(t)=q(t)+\frac{C_P(\mu)}{Z\,t}v(t) $. We want to show that $f\leq 0$. Compute: $$f'(t)=\frac{1}{Z\,h(t)}\exp\left(-\lambda_\mu \int_0^t\frac{u\,du}{h(u)}\right)-\frac{C_P(\mu)}{Zt^2}v(t)-\lambda_\mu \frac{tC_P(\mu)}{Z\,h(t)t}v(t)=-\frac{v(t)}{Z}\frac{C_P(\mu)}{t^2} <0.$$ Hence $f\leq \underset{t\rightarrow a}{\lim}f(t)\leq0$.

\end{enumerate}
\end{dem}

\section{Exact and approximate integration by parts formulas}\label{sectipp}

Let us recall a classical result in $\Gamma$-calculus which we will often use in the sequel. For all $f,g\in H^1(\mu)$, the following integration by parts formula holds (see formula ($\ref{defH1}$) for the definition of $H^1(\mu)$).
\begin{equation}\label{ippclassique}
\int\Gamma(f,g)d\mu = - \int fLg\,d\mu
\end{equation}
Indeed, since $\mu$ is the invariant measure of the process generated by $L$, $\int L\phi\,d\mu = 0$ for all functions $\phi \in \mathcal{D}$. Hence, integrating the definition of the carré du champ operator $\Gamma(f,g):=\frac{1}{2}\left(L(fg)-fLg-gLf\right) $ and using the reversibility of $\mu$, one gets the result.
Let us now state an integration by parts formula for the quotient process.
\begin{prop}\label{ippmustar}
For all $C^1(I)$ functions $\psi:I\rightarrow \mathbb{R}$ such that $\psi\circ f_0\in H^1(\mu)$, it holds that: 
\begin{equation}\label{eqippmustar}
\int_{\mathbb{R}}x\psi(x)d\mu^*(x) = C_{P}(\mu) \int_\mathbb{R}h(x)\psi'(x)d\mu^*(x)
\end{equation}
\end{prop}

\begin{dem}
As a consequence of the integration by part formula ($\ref{ippclassique}$) for the initial diffusion process, one gets $$ \int \Gamma(f_0,g)d\mu = -\int g Lf_0 \, d\mu = \frac{1}{C_P(\mu)}\int g f_0\, d\mu.$$
We then use this equality with $g:=\phi\,\circ f_0$ and use the diffusion property of $\Gamma$ and the definition of $h$.
\end{dem}
Let us state now an extension of the approximate integration by part formula which extend a Lemma from Courtade and Fathi \cite[Lemma 2.3]{CF}. Let us recall that $$H^1(\nu):=\{f\in L^2(\nu)\cap\mathcal{D}\ |\ \int f\, d\nu = 0 \, \mathrm{and}\, \int\Gamma (f) d\nu<\infty\} .$$

\begin{thm}\label{ippapprox}
The following inequality holds for all $g\in H^1(\nu)$: 
\begin{equation}\label{eqippapprox}
\left| \int f_0 g\ d\nu - C_p(\nu)\int \Gamma(f_0,g)d\nu \right| \leq \left( C_P(\nu) - C_P(\mu)\right)^\frac{1}{2}\left(\frac{C_P(\nu)}{C_P(\mu)}\right)^{\frac{1}{2}}\left(\int\Gamma (g)\ d\nu\right)^\frac{1}{2} 
\end{equation}
\end{thm}

\begin{dem}
Let $t\in\mathbb{R}$ and $g\in H^1(\nu)$. Let us apply the Poincaré inequality for $\nu$ to $\alpha:= f_0 + tg$. Computing that $$ \mathrm{Var}_{\nu}(\alpha) = \int (f_0+tg)^2 \,d\nu = \underset{=1}{\underbrace{\int f_0^2d\nu}} + 2t\int f_0 g\,d\nu + t^2\int g^2d\nu,$$ and that $$\int\Gamma(\alpha)d\nu = \underset{\leq C_P(\mu)^{-1}}{\underbrace{\int \Gamma(f_0)d\nu}}  + 2t\int\Gamma(f_0,g)d\nu +t^2\int\Gamma(g)d\nu, $$ we get that $$ 1 + 2t\int f_0 g\,d\nu + t^2\int g^2d\nu \leq C_P(\nu)\left( C_P(\mu)^{-1} + 2t\int\Gamma(f_0,g)d\nu +t^2\int\Gamma(g)d\nu  \right). $$ 

Now, writting $C_P(\nu) = C_P(\nu)-C_P(\mu)+C_P(\mu) =: \Delta C_P + C_P(\mu)$, we obtain : $$ 2t\int f_0 g\,d\nu + t^2\int g^2d\nu \leq \frac{\Delta C_P}{C_P(\mu)}+  2C_P(\nu)t\int\Gamma(f_0,g)d\nu +t^2C_P(\nu)\int\Gamma(g)d\nu.  $$ Hence, for all $t\in \mathbb{R}$, $$ \frac{\Delta C_P}{C_P(\mu)}+  2t\left(C_P(\nu)\int\Gamma(f_0,g)d\nu -\int f_0 g\,d\nu\right) +t^2C_P(\nu)\int\Gamma(g)d\nu \geq 0. $$  This polynomial of degree $2$ takes only non-negative values, hence its discriminant is non-positive: $$b^2-4ac=4\left(C_P(\nu)\int\Gamma(f_0,g)d\nu-\int f_0 g\,d\nu\right)^2 - 4C_P(\nu)\int\Gamma(g)d\nu  \frac{\Delta C_P}{C_P(\mu)} \leq 0,$$ which yields the result.
\end{dem}

Using Theorem $\ref{ippapprox}$ with $g:=\psi\circ f_0$ and dividing by $C_P(\nu)$, we can state now the approximate integration by part formula for the quotient process.

\begin{corol}\label{ippnustar}
For all $C^1(I)$ functions $\psi:I\rightarrow \mathbb{R}$ such that $\psi\circ f_0\in H^1(\mu)$, $$ \left| \int \left( h(x)\psi'(x) - \frac{x}{C_P(\nu)}\psi(x)\right) d\nu^* \right| \leq \left( \frac{ C_P(\nu) - C_P(\mu)}{C_P(\nu)\,C_p(\mu)} \right)^\frac{1}{2}\left(\int h(x)\psi'(x)^2 d\nu^*\right)^\frac{1}{2}.  $$
\end{corol}

The original point of view of Courtade and Fathi is the following \cite[Section 1.3]{CF}. The Poincaré constant is considered as the minimizer of the energy $\int \Gamma(f)\,d\mu$ for $f$ with variance $1$ in a large enough functional set. One can then write the Euler-Lagrange equation. The heuristic is that if another measure almost satisfies this equation, its Poincaré constant would not be so far from the one of $\mu$. This is made rigorous with Stein's method.
Here the integration by part formula ($\ref{eqippmustar}$) plays the role of the Euler-Lagrange equation.

\section{Implementing Stein's method}\label{steinmethod}

The original idea of Stein's method is to control some distance between two probability measures by bounding the solution of some equation called Stein equation. For more details, see the survey \cite{fundamStein}.
We start with the approximated integration by parts formula given in Corollary $\ref{ippnustar}$ above.
The goal is to get stability inequalities of the form $$C_P(\nu) \geq C_P(\mu) + \alpha\, d(\mu^*,\nu^*)^2, $$
where $\alpha>0$ is a multiplicative constant and $d$ is a metric on the space of probability measures defined by
$$ d(\mu,\nu) := \underset{f\in \mathcal{F}}{\sup}\left| \int fd\mu - \int f d\nu \right|  $$
for some set of functions $\mathcal{F}$. Metrics of this form include $1$-Wasserstein, Kolmogorov and total variation distances.
Denoting $\mu^*(f):=\int f\,d\mu^*$, this goal would be achieved if one could solve the equation $$ S_\nu(\psi)(x):= h(x)\psi'(x)-\frac{x}{C_P(\nu)}\psi(x)=f(x)-\mu^* (f) \ , \ x\in I,  $$ for all $f\in\mathcal{F}$, and bound the term $\int h(x)\psi'(x)^2 d\nu^*$ independently of $f\in\mathcal{F}$.
The actual Stein equation is $$ S_\mu(\psi)(x):=h(x)\psi'(x)-\frac{x}{C_P(\mu)}\psi(x)=f(x)-\mu^* (f) \ , \ x\in \overset{\circ}{I} . $$ Since $S_\mu(\psi')=L^*(\psi)$, for any solution $\phi$ of the Poisson equation $L\phi=f-\mu^*(f)$, $\phi'$ solves the Stein equation. Hence one can study the Stein equation via the probabilistic analysis of the Poisson equation. Observing that $$ S_\mu \psi = S_\nu \psi + \left( \frac{1}{C_P(\nu)}-\frac{1}{C_P(\mu)} \right)\, x\,\psi, $$ Corollary $\ref{ippnustar}$ gives
$$\left| \int S_\mu(\psi)\, d\nu^* \right| \leq \left( \frac{ C_P(\nu) - C_P(\mu)}{C_P(\nu)\,C_p(\mu)} \right)^\frac{1}{2}\left(\int h(x)\psi'(x)^2 d\nu^*\right)^\frac{1}{2} + \left( \frac{ C_P(\nu) - C_P(\mu)}{C_P(\nu)\,C_p(\mu)} \right)\left|\int x\, \psi\, d\nu^*\right|.  $$
Now, using that $\nu^*$ is centered, the Cauchy-Schwarz inequality and the Poincaré inequality for $\nu^*$, we obtain, denoting $\mathrm{Var}_{\nu^*}(\psi):= \int (\psi-\int\psi d\nu^*)^2 d\nu^*$,
\begin{align*}
\left|\int x\, \psi\, d\nu^*\right| &= \left|\int x\, (\psi-\int\psi d\nu^*)\, d\nu^*\right| \\
&\leq \left( \mathrm{Var}_{\nu^*}(\psi) \right)^{\frac{1}{2}} (\int x^2\,d\nu^*)^{\frac{1}{2}} \\
 &= \left( \mathrm{Var}_{\nu^*}(\psi) \right)^{\frac{1}{2}} (\int f_0^2 d\nu)^{\frac{1}{2}} \\
 &= \left( \mathrm{Var}_{\nu^*}(\psi) \right)^{\frac{1}{2}}\\
 &\leq \sqrt{C_P(\nu)}\left(\int \Gamma^*(\psi)\, d\nu^*\right)^\frac{1}{2}.
\end{align*}
Setting, $$\delta:= \frac{ C_P(\nu) - C_P(\mu)}{C_P(\nu)\,C_P(\mu)}$$ we have proved that
\begin{eqnarray}\label{ippaproxresume}
 \left| \int S_\mu(\psi)\, d\nu^* \right| \leq \left(\sqrt{\delta} + \sqrt{C_P(\nu)}\,\delta \right) \left(\int \Gamma^*(\psi)\, d\nu^*\right)^\frac{1}{2} . 
\end{eqnarray}
So our goal would be achieved if for all $f\in \mathcal{F}$, we can solve the Stein equation $S_\mu(\psi)=f-\mu^*(f)$, and get a bound for $\int \Gamma^*(\psi)\,d\nu^*=\int_a^b h(x)\psi'(x)^2 d\nu^*$ independently of $f$. Note that it would be enough to get a bound for $||\sqrt{h}\psi'||_\infty$, or for $||\psi'||_\infty$ since $\int h(x)\psi'(x)^2 d\nu^*\leq ||\psi'||_\infty^2\int h\,d\nu^* = ||\psi'||_\infty^2 \int \Gamma(f_0)\,d\nu \leq \frac{||\psi'||_\infty^2}{C_P(\mu)} $, where $||\cdot ||_\infty $ denotes the supremum norm .

\subsection{Stein solution}\label{steinsolsection}

Let $f\in L^1(\mu^*) $. We consider the Stein equation 
\begin{equation}\label{steineq}
h(x)\psi'(x)-\frac{x}{C_P(\mu)}\psi(x)=f(x)-\mu^* (f) \ , \ x\in I 
\end{equation} 
In the sequel, we will call $f$ the target function. This equation is a first order linear ODE, which can be explicitly solved.
\begin{lemm}\label{solstein}
The function 
\begin{equation}\label{steinsolgauche}
\psi(x):=\exp\left(  \int_{0}^{x}\frac{u\,du}{C_P(\mu) h(u)}\right)\int_{a}^{x}\left(f(y)-\mu^* (f)\right) \frac{1}{ h(y)}\exp\left( - \int_{0}^{y}\frac{u\,du}{C_P(\mu) h(u)}\right)  dy 
\end{equation} is a solution of ($\ref{steineq}$). Moreover, it also can be written as 
\begin{equation}\label{steinsoldroite}
\psi(x)= - \exp\left(  \int_{0}^{x}\frac{u\,du}{C_P(\mu) h(u)}\right)\int_{x}^{b}\left(f(y)-\mu^* (f)\right) \frac{1}{ h(y)}\exp\left( - \int_{0}^{y}\frac{u\,du}{C_P(\mu) h(u)}\right) dy 
\end{equation} 
\end{lemm}

\begin{dem}
The solution ($\ref{steinsolgauche}$) is obtained by classical tools for first order linear differential equations. The second formula is obtained from Proposition $\ref{densitymustar}$.
\end{dem}
Let us emphasize that this approach is only possible in dimension $1$.

\subsection{Boundedness of solutions for bounded target functions}

Under the current framework one can give the following bounds for the solution of the Stein equation.

\begin{prop}\label{solutioncontrol2}
Let $f$ be a bounded measurable function on $I$, and let $\psi$ be the solution of $(\ref{steineq})$ given by Lemma $\ref{solstein}$. Then 
\begin{equation*}
||\psi||_\infty \leq Z\,\max (q(0), 1-q(0))\,||f-\mu^* (f)||_\infty \leq Z\, ||f-\mu^*(f)||_{\infty},
\end{equation*} with $Z$ the normalization constant (see Proposition $\ref{densitymustar}$) and $q$ the cumulative distribution of $\mu^*$. Moreover,
\begin{equation*}
||x\psi(x)||_\infty \leq C_P(\mu)\,||f-\mu^*(f)||_{\infty}.
\end{equation*}
\end{prop}

\begin{dem}
Let us recall that $v(t) = \exp\left( -\lambda_\mu\int_0^t \frac{u\,du}{h(u)} \right)$.
If $x\leq 0$ we use the writting ($\ref{steinsolgauche}$) of $\psi$, and get: $$ |\psi(x)|\leq ||f-\mu^* (f)||_\infty v(x)^{-1}\,Z\,q(x). $$ and by Lemma $\ref{tailcontrol}$: $$ \underset{x\in I_-}{\sup}|\psi(x)| \leq Z\, q(0)\,||f-\mu^* (f)||_\infty. $$ If $x>0$, we use the writting ($\ref{steinsoldroite}$) of $\psi$ in Lemma $\ref{solstein}$ together with Lemma $\ref{tailcontrol}$: $$\underset{x\in I_+}{\sup}|\psi(x)| \leq Z\, (1-q(0))\,||f-\mu^* (f)||_\infty. $$

To prove the second inequation, we split off again the case $x>0$ and $x<0$ using in each case the appropriate representation of $\psi$ in Lemma $\ref{solstein}$ and we conclude with Lemma $\ref{tailcontrol}$.
If $x\in I_{-}$, then $|\psi(x)|\leq v(x)^{-1}||f-\mu^*(f)||_\infty \, Z\, q(x) $ but $q(x)\leq \frac{C_P(\mu)}{Z\,|x|}v(x)$ so $|x\psi(x)|\leq C_P(\mu)\,||f-\mu^*(f)||_{\infty}$.
If $x\in I_+$, then $|\psi(x)|\leq v(x)^{-1}||f-\mu^*(f)||_\infty \, Z\, (1-q(x)) $ but $(1-q(x))\leq \frac{C_P(\mu)}{Z\,x}v(x)$ so $|x\psi(x)|\leq C_P(\mu)\,||f-\mu^*(f)||_{\infty}$.
\end{dem}

\section{Stability in total variation under a uniform ellipticity condition}\label{ellipt}

Let us recall the definition of the total variation distance.
\begin{defi}
The total variation distance between two measures $\alpha, \beta $ on $I$ is $$d_{TV}(\alpha,\beta) = \underset{A\subset I}{\sup} \left|\int \textbf{1}_A d\alpha - \int \textbf{1}_A d\beta \right|, $$ where the supremum is running over all measurable subsets of $I$.
\end{defi}

We assume the following uniform ellipticity condition on $h$:
\begin{equation}\label{hminore}
\kappa:= \underset{x\in\mathbb{R}}{\inf} h(x) >0
\end{equation}
Under this ellipticity condition, one can get better estimates on the solutions of Stein equation. The vocabulary of ellipticity comes from PDE theory, see \cite[Chapter 6]{pde}. This condition extends the Gaussian case discussed at the beginning of Section $\ref{quotienprocess}$. We can now state the main result of this section.

\begin{thm}\label{stabilityhminore}
Let $(L,\mu)$ be a Markov diffusion on a Riemannian manifold $M$ satisfying a Poincaré inequality with sharp constant $C_P(\mu)$. Suppose that Assumptions $\ref{spectralgapattained}$, $\ref{gammaf0measuable}$ and $\ref{criticalmaxima}$ are verified. If $\kappa:= \underset{x\in\mathbb{R}}{\inf} h(x) >0$, then for all measures $\nu$ on $M$ normalized as in ($\ref{hypnu}$) and satisfying a Poincaré inequality with sharp constant $C_P(\nu)$,
\begin{equation}\label{stabilitytv}
d_{TV}(\mu^*,\nu^*) \leq \frac{4}{\sqrt{\kappa}}\left(\sqrt{\delta} + \sqrt{C_P(\nu)}\,\delta \right),
\end{equation}
where $\mu^*$ (resp. $\nu^*$) is the pushforward of $\mu$ (resp. $\nu$) by the first eigenfunction $f_0$, and $\delta:= \frac{ C_P(\nu) - C_P(\mu)}{C_P(\nu)\,C_p(\mu)}$.
\end{thm}

In order to prove Theorem $\ref{stabilityhminore}$, we need the following bounds on the solution of the Stein equation.

\begin{prop}\label{gammasolutioncontrol}
Let $f$ be a bounded function on $\mathbb{R}$, and let $\psi$ be the solution of the Stein equation given by Lemma $\ref{solstein}$. Then 
$$ ||\psi'||_\infty \leq \frac{2}{\kappa }\,||f-\mu^* (f)||_\infty, $$ and
$$ ||h\psi'^2||_\infty \leq \frac{4}{\kappa }\,||f-\mu^* (f)||_{\infty}^2. $$
\end{prop}

\begin{dem}
By direct calculation: $$\psi'(x) = \frac{1}{ h(x)}\left( f - \mu^* (f) + \frac{x}{C_P(\mu)}\psi(x)\right).$$ The lower bound on $h$ and Proposition $\ref{solutioncontrol2}$ yield the first claim.
Similarly, $$h(x)\psi'(x)^2 = \frac{1}{ h(x)}\left( f - \mu^* (f) + \frac{x}{C_P}\psi(x)\right)^2.$$
So from the elementary inequality $(a+b)^2\leq 2a^2+2b^2$, and the same reasoning as above, we get the result.
\end{dem}

\begin{demof}{\textbf{Theorem} $\textbf{\ref{stabilityhminore}}$}
Let $f:I\rightarrow \mathbb{R}$ be an indicator function.
We solve the Stein equation with Lemma $\ref{solstein}$ and denote $\psi$ the chosen solution. Then ($\ref{ippaproxresume}$) and Proposition $\ref{gammasolutioncontrol}$ give: $$\left| \int S_\mu(\psi)\, d\nu^* \right| \leq \left(\sqrt{\delta} + \sqrt{C_P(\nu)}\,\delta \right) \left(\int \Gamma^*(\psi)\, d\nu^*\right)^\frac{1}{2} \leq \left(\sqrt{\delta} + \sqrt{C_P(\nu)}\,\delta \right)\frac{2}{\sqrt{\kappa}}||f-\mu^* (f)||_{\infty}.$$ But: $$ \left| \int S_\mu(\psi)\, d\nu^* \right| = \left| \int\left( f - \mu^*(f) \right) d\nu^* \right| = \left|\nu^*(f) - \mu^*(f) \right|. $$
Hence for any indicator function $f$, 
$$ \left|\nu^*(f) - \mu^*(f) \right| \leq \frac{4}{\sqrt{\kappa}}\left(\sqrt{\delta} + \sqrt{C_P(\nu)}\,\delta \right). $$
Taking the supremum over all such $f$ concludes the proof.
\end{demof}
Let us point out that ($\ref{stabilitytv}$) is of the form ($\ref{stabdef}$) with $\phi(x,y)=\frac{4}{\sqrt{\kappa}}\left( \sqrt{\frac{y-x}{xy}}+ \frac{y-x}{x\sqrt{y}}\right)$.\\

When specialized to the Gaussian case, Theorem $\ref{stabilityhminore}$ gives that if $\nu$ is a normalized probability measure on $\mathbb{R}$ and $C_P(\nu)-1\leq 1$, then:
$$C_P(\nu)\geq 1+\frac{1}{16\left(1+\sqrt{2}\right)} d_{TV}(\gamma,\nu)^2 . $$ Our constant $16\left(1+\sqrt{2}\right)$ is worse than the constant $9$ from \cite{U}.

\section{Stability in $W_1$ distance}\label{explicit}

Let us recall the Kantorovitch-Rubinstein formula, which we will take as the definition of the $W_1$ distance.
\begin{defi} (see \cite[Chapter 5]{villani}  )\label{w1int}
The $1$-Wasserstein distance between two measures $\alpha, \beta $ on $I$ is $$W_1(\alpha,\beta) = \underset{\underset{1-\mathrm{Lipschitz}}{f:I\rightarrow \mathbb{R} }}{\sup} \left|\int fd\alpha - \int f d\beta \right|. $$ 
\end{defi} 
We can now state the main theorem of this article.
\begin{thm}\label{thmcentr}
Let $(L,\mu)$ be a Markov diffusion on a Riemannian manifold $M$ satisfying a Poincaré inequality with sharp constant $C_P(\mu)$. Suppose that Assumptions $\ref{spectralgapattained}$, $\ref{gammaf0measuable}$ and $\ref{criticalmaxima}$ are verified. Let
$$
C_h := \underset{x\in I}{\sup} \left[ \sqrt{\Gamma^*(a_1)(x)}\int_{a}^x q(t)\,dt +  \sqrt{\Gamma^*(a_2)(x)}\int_{x}^b (1-q(t))\,dt \right] ,
$$
where $a_1$ and $a_2$ are defined in ($\ref{defa1}$) and ($\ref{defa2}$). Then for all measures $\nu$ on $M$ normalized as in ($\ref{hypnu}$) and satisfying a Poincaré inequality with sharp constant $C_P(\nu)$, it holds:
$$W_1(\mu^*,\nu^*) \leq C_h\left(\sqrt{\delta} + \sqrt{C_P(\nu)}\,\delta \right) \leq C_h\left( \frac{\sqrt{C_P(\nu)-C_P(\mu)}}{C_P(\mu)}+ \frac{C_P(\nu)-C_P(\mu)}{C_P(\mu)^{3/2}}\right), $$
where $\mu^*$ (resp. $\nu^*$) is the pushforward of $\mu$ (resp. $\nu$) by $f_0$, $W_1$ is the $1$-Wasserstein distance (see Definition $\ref{w1int}$), and $\delta:= \frac{ C_P(\nu) - C_P(\mu)}{C_P(\nu)\,C_p(\mu)}$. 
\end{thm}
\begin{dem}
The proof is the same as in Theorem $\ref{stabilityhminore}$ but we replace the indicator functions by $1$-Lipschitz functions, and we use Proposition $\ref{gammasteinbound}$ instead of Proposition $\ref{gammasolutioncontrol}$.
\end{dem}
All the sequel of this section is devoted to prove Proposition $\ref{gammasteinbound}$.
In \cite[Appendix p. 37]{CGS}, Chen, Goldstein and Shao bound the derivatives of the solution of the Stein equation in the case where $\mu$ is the standard Gaussian measure.
Here, we extend their approach to a broader class of measures. Let $f:I\rightarrow \mathbb{R}$ be absolutely continuous, and let $\psi$ be the solution of ($\ref{steineq}$) given by Lemma $\ref{solstein}$. We are looking for inequalities of the form $$||h\psi'^2||_\infty \leq C\, ||f'||_\infty, $$ where $C>0$ is a constant which does not depend on $f$.
Recall that we denote by $q(t)=\int_a^t d\mu^*$ the cumulative distribution function of the measure $\mu^*$, and $v(t)=\exp\left( -\lambda_\mu \int_{0}^{x}\frac{u\,du}{h(u)}\right)$ (see ($\ref{defdev}$)).

\begin{prop}\label{solreecri}
Let $f:I\rightarrow \mathbb{R}$ be in $ C^1(I)\cap L^1(\mu^*)$. 
The associated solution ($\ref{steinsoldroite}$) of Lemma $\ref{solstein}$ can be written as: 
\begin{equation}\label{nvpsi}
 \psi(x)=-Z \frac{1-q(x)}{v(x)}\int_{t=a}^x f'(t)q(t)\,dt \, -\, Z \frac{q(x)}{v(x)}\int_{t=x}^b  f'(t)(1-q(t))\,dt
\end{equation}
\end{prop}

\begin{dem}
First, write 
\begin{align*}
f(x)-\mu^*(f) &= \int_a^x (f(x)-f(y))\,d\mu^*(y) - \int_b^x (f(x)-f(y))\,d\mu^*(y) \\
&= \int_a^x \int_y^xf'(t)\,dt\,d\mu^*(y) - \int_b^x \int_y^xf'(t)\,dt\,d\mu^*(y) \\
&= \int_a^xf'(t)\int_a^td\mu^*(y)\,dt - \int_x^bf'(t)\int_t^bd\mu^*(y)\,dt
\end{align*}
\begin{align}\label{formulainterme}
&= \int_a^xf'(t)q(t)\,dt - \int_x^bf'(t)(1-q(t))\,dt.
\end{align}
It follows that 
\begin{align*}
\psi(x) &= \frac{1}{v(x)}\int_a^x(f(y)-\mu^*(f))Z\,d\mu^*(y)\\
&= \frac{1}{v(x)}\int_a^x \left(\int_a^yf'(t)q(t)\,dt - \int_y^bf'(t)(1-q(t))\,dt\right)Z\,d\mu^*(y).
\end{align*}
Now, on the one hand: $$ \int_a^x \int_{a}^y f'(t) q(t) \,dt\,d\mu^*(y)  =  \int_{a}^x \int_{t}^x f'(t)q(t)\,d\mu^*(y)\,dt   
 =\int_{a}^x  f'(t) q(t)\left(q(x)-q(t) \right)dt, $$ on the other hand: 
\begin{align*}
\int_a^x \int_{y}^b f'(t)(1-q(t))\,dt \,d\mu^*(y)
 &=\int_{x}^b \int_{a}^x f'(t)(1-q(t)) \,d\mu^*(y)\,dt 
 + \int_{a}^x \int_{a}^t f'(t)(1-q(t))\,d\mu^*(y)\,dt\\ 
 &=\int_{x}^b  f'(t)(1-q(t))q(x) \,dt
 + \int_{a}^x  f'(t)(1-q(t)) q(t)\,dt.
\end{align*}  
From which one gets ($\ref{nvpsi}$).
\end{dem}
By definition $\psi$ satisfies $$h(x)\psi'(x) = \left( f(x) - \mu^*(f) + \frac{x}{C_P}\,\psi(x) \right). $$
Hence, using Proposition $\ref{solreecri}$ and Formula ($\ref{formulainterme}$), we are now able to give a new formula for the derivative of the solution.

\begin{corol}\label{derivsolreecri}
We have:
\begin{align*}
h(x)\psi'(x) &= \left(1-Z (1-q(x))\frac{x}{C_P(\mu)}\exp\left( \int_{0}^{x}\frac{u\,du}{C_P(\mu) h(u)}\right)\right)\int_{a}^x  f'(t) q(t) \,dt \\
 &- \left(1+Z q(x)\frac{x}{C_P(\mu)}\exp\left( \int_{0}^{x}\frac{u\,du}{C_P(\mu) h(u)}\right)\right) \int_{x}^b  f'(t)(1-q(t))\,dt.
\end{align*}
\end{corol} 
In order to simplify the expression above, one can see that 
\begin{equation}\label{gammaa1}
\left| 1-Z (1-q(x))\frac{x}{C_P(\mu)}\exp\left( \int_{0}^{x}\frac{u\,du}{C_P(\mu) h(u)}\right)\right|\frac{1}{\sqrt{h(x)}} = \sqrt{\Gamma^*(a_1)(x)}
\end{equation} and
\begin{equation}
\left|1+Z q(x)\frac{x}{C_P(\mu)}\exp\left( \int_{0}^{x}\frac{u\,du}{C_P(\mu) h(u)}\right)\right|\frac{1}{\sqrt{h(x)}} =  \sqrt{\Gamma^*(a_2)(x)}
\end{equation}
where 
\begin{equation}\label{defa1}
a_1(x):= Z(1-q(x))\exp\left(\int_0^x \frac{u\,du}{C_P(\mu)h(u)}\right)
\end{equation}

\begin{equation}\label{defa2}
a_2(x):= Zq(x)\exp\left(\int_0^x \frac{u\,du}{C_P(\mu)h(u)}\right)
\end{equation}
The following is then immediate.

\begin{prop}\label{gammasteinbound}
Let $f:I\rightarrow \mathbb{R}$ be in $ C^1(I)\cap L^1(\mu^*)$, and let $\psi$ the associated solution ($\ref{steinsoldroite}$). The following bound holds:
$$||\sqrt{h}\psi'||_\infty \leq C_h\, ||f'||_\infty, $$
where 
\begin{equation}\label{defdeCh}
C_h := \underset{x\in I}{\sup} \left[ \sqrt{\Gamma^*(a_1)(x)}\int_{a}^x q(t)\,dt +  \sqrt{\Gamma^*(a_2)(x)}\int_{x}^b (1-q(t))\,dt \right]
\end{equation}
\end{prop}

\subsection{Finitness of $C_h$}\label{finitch}
If $C_h=\infty$, the result of Theorem $\ref{thmcentr}$ degenerates. 
That is why in this section, we give explicit conditions on $h$ that ensure the finitness of $C_h$.

\begin{thm}\label{thmcentralexploitable}
Let $(L,\mu)$ be a Markov diffusion on a Riemannian manifold $M$ satisfying a Poincaré inequality with sharp constant $C_P(\mu)$. Suppose that Assumptions $\ref{spectralgapattained}$, $\ref{gammaf0measuable}$ and $\ref{criticalmaxima}$ are verified.  Recall that we denote $a=\inf f_0(M)<0$ and $b=\sup f_0(M)>0$. Assume that one of these two conditions is verified at $a$: 
\begin{itemize}
\item either $a=-\infty$ and $c_1 |t|^{2\alpha-2}\leq h(t)\leq c_2 |t|^\alpha$ for $t\rightarrow -\infty$ with $\alpha\leq 2 $ and $c_1,\,c_2>0$,
\item or $a>-\infty$ and $c_1 (t-a)^2\leq h(t)\leq c_2(t-a)$ for $t\rightarrow a^+$ with $c_1,\,c_2>0$,
\end{itemize}
and one of these two conditions is satisfied at $b$:
\begin{itemize}
\item either $b=+\infty$ and $c_1 t^{2\alpha-2}\leq h(t)\leq c_2 t^\alpha$ for $t\rightarrow +\infty$ with $\alpha\leq 2 $ and $c_1,\,c_2>0$,
\item or $b<+\infty$ and $c_1 (b-t)^2\leq h(t)\leq c_2(b-t)$ for $t\rightarrow b^-$ with $c_1,\,c_2>0$.
\end{itemize}
Then for all measures $\nu$ on $M$ normalized as in ($\ref{hypnu}$) and satisfying a Poincaré inequality with sharp constant $C_P(\nu)$, it holds:
$$W_1(\mu^*,\nu^*) \leq C_h\left(\sqrt{\delta} + \sqrt{C_P(\nu)}\,\delta \right), $$
where $C_h>0$ is a constant, $\mu^*$ (resp. $\nu^*$) is the pushforward of $\mu$ (resp. $\nu$) by $f_0$, $W_1$ is the $1$-Wasserstein distance (see Definition $\ref{w1int}$), and $\delta:= \frac{ C_P(\nu) - C_P(\mu)}{C_P(\nu)\,C_p(\mu)}$. 
\end{thm}

In order to prove Theorem $\ref{thmcentralexploitable}$, we will use the following results.

\begin{thm}\label{compath}
Assume that:
\begin{enumerate}
\item There exists a non-negative function $g_1$, defined on a neighborhoood of $a$, satisfying 
\begin{equation}\label{bornitude}
\frac{x(g_1(x)-1)}{\sqrt{h(x)}}=O(1),\quad x\rightarrow a^+
\end{equation}  and  
\begin{equation}\label{g_1} q(x)\geq-\frac{C_P(\mu)}{Z}v(x)\frac{g_1(x)}{x}, \quad x\rightarrow a^+
\end{equation}
\item There exists a non-negative function $g_2$, defined on a neighborhoood of $b$, satisfying $$\frac{x(g_2(x)-1)}{\sqrt{h(x)}}=O(1),\quad x\rightarrow b^-$$ and
\begin{equation}\label{g_2}
1-q(x)\geq \frac{C_P(\mu)}{Z}v(x)\frac{g_2(x)}{x}, \quad x\rightarrow b^-
\end{equation}
\end{enumerate}
Then the constant $C_h$ defined in Theorem $\ref{thmcentr}$ is finite.
\end{thm}

\begin{dem}
We only study boundedness of $\sqrt{\Gamma^*(a_1)(x)}\int_{-\infty}^x q(t)\, dt $ at $a$ and $b$ and conclude by continuity. The case of $\sqrt{\Gamma^*(a_2)(x)}\int_{x}^{+\infty} (1-q(t))\,dt $ is similar. We only deal with the case of $a$, the one of $b$ is similar. Since $v'(x)=-\lambda_\mu h(x)v(x)$ with an integration by parts, we get 
$$\int_{a}^x q(t)\,dt =\left[tq(t)\right]_{a}^x - \int_{a}^x \frac{1}{Z}t\frac{v(t)}{h(t)}\,dt \leq  xq(x)+\frac{C_P(\mu)}{Z}v(x). $$ The last inequality comes from Lemma $\ref{tailcontrol}$ which gives $\underset{t\rightarrow a}{\lim}\left(-tq(t)-\frac{C_P(\mu)}{Z}v(t)\right)\leq 0 $. Moreover, it is obvious that $1-Z \lambda_\mu(1-q(x)) x\exp\left( \lambda_\mu \int_{0}^{x}\frac{u\,du}{ h(u)}\right) \geq 0$, for $x<0$. Hence for $x<0$, using ($\ref{gammaa1}$),
\begin{align*}
\sqrt{\Gamma^*(a_1)(x)}\int_{a}^x q(t)\, dt &= \frac{1}{\sqrt{h(x)}}\left(1-Z \lambda_\mu(1-q(x)) x\exp\left( \lambda_\mu \int_{0}^{x}\frac{u\,du}{ h(u)}\right) \right) \int_{a}^x q(t)\, dt   \\
&\leq \frac{1}{\sqrt{h(x)}}\left(xq(x)+\frac{C_P(\mu)}{Z}v(x)\right) - Z\lambda_\mu\frac{1-q(x)}{\sqrt{h(x)}}\frac{x^2}{v(x)}q(x)-\frac{1-q(x)}{\sqrt{h(x)}}x. 
\end{align*} 
On the one hand, using first the assumption on $g_1$, and next Lemma $\ref{tailcontrol}$, one gets that for $x\rightarrow a^+$, $$ \frac{1}{\sqrt{h(x)}}\left(xq(x)+\frac{C_P(\mu)}{Z}v(x)\right) \leq -\frac{C_P(\mu)}{Z\sqrt{h(x)}}v(x)g_1(x)+\frac{C_P(\mu)}{Z\sqrt{h(x)}}v(x) \leq \frac{(g_1(x)-1)xq(x)}{\sqrt{h(x)}} \underset{t\rightarrow a}{\rightarrow} 0. $$
On the other hand, the existence of $g_1$ gives $$ -  - Z\lambda_\mu\frac{1-q(x)}{\sqrt{h(x)}}\frac{x^2}{v(x)}q(x)-\frac{1-q(x)}{\sqrt{h(x)}}x\leq x\frac{g_1(x)}{\sqrt{h(x)}}(1-q(x))-\frac{1-q(x)}{\sqrt{h(x)}}x = \frac{x(g_1(x)-1)}{\sqrt{h(x)}}(1-q(x)), $$ which is bounded by assumption.\\	
The boundedness at $b$ is obtained in the same way. By Lemma $\ref{tailcontrol}$, one gets that $1-Z\lambda_\mu (1-q(x))x\exp\left( \lambda_\mu \int_{0}^{x}\frac{u\,du}{ h(u)}\right) \geq 0$, for $x>0$. By the assumption on $g_2$, one gets $$1-Z \lambda_\mu(1-q(x)) x\exp\left( \lambda_\mu \int_{0}^{x}\frac{u\,du}{ h(u)}\right) \leq 1-g_2(x).$$ Hence using the same integration by parts as before, $$\sqrt{\Gamma^*(a_1)(x)}\int_{a}^x q(t)\, dt \leq \frac{1}{\sqrt{h(x)}}\left(xq(x)+\frac{C_P(\mu)}{Z}v(x)\right)(1-g_2(x)).$$ Now $\frac{1-g_2(x)}{\sqrt{h(x)}}v(x)=\frac{x(1-g_2(x))}{\sqrt{h(x)}}\frac{v(x)}{x}$ goes to zero when $x$ goes to $b$ since ($\ref{g_2}$) gives $\frac{v(x)}{x}\underset{x\rightarrow b}{\rightarrow} 0$. This concludes the proof.
\end{dem}
The standard Gaussian satisfies the requirements of Theorem $\ref{compath}$ with  $g_1(x)=g_2(x)=\frac{x^2}{1+x^2}$.\\
Let us now give sufficient conditions on $h$ to ensure that such $g_1,g_2$ functions exist. We prove the case of $g_1$ at $a$, the case of $g_2$ at $b$ is similar.

\begin{lemm}\label{lemmcondipourg1}
Let $g_1$ be a bounded non-negative $C^1$ function defined on a neighborhood of $a$ such that for $t\rightarrow a^+ $, $g_1(t) >tg_1'(t)$ and 
\begin{equation}\label{condipourg1}
h(t)\leq \lambda_\mu \frac{t^2(1-g_1(t))}{g_1(t)-tg_1'(t)},\quad t\rightarrow a^+
\end{equation}
If $\underset{t\rightarrow a}{\lim}\,v(t)=0 $, then $$q(t)\geq-\frac{C_P}{Z}v(t)\frac{g_1(t)}{t},\quad t\rightarrow a^+ . $$
\end{lemm}
\begin{dem}
Let $f(t):=q(t)+\frac{C_P}{Z}\frac{g_1(t)}{t}v(t)$. Compute $$f'(t)=\frac{v(t)}{Z}\left( \frac{C_P\,g_1'(t)}{t}- \frac{C_P\,g_1(t)}{t^2}+\frac{1}{h(t)}-\frac{g_1(t)}{h(t)} \right) $$ The assumption implies then that $f'\geq 0$. Hence $f(t)\geq \underset{x \rightarrow a}{\lim} f(x) = 0$.
\end{dem}

\begin{lemm}\label{lemmcondipourg2}
Let $g_2$ be a bounded non-negative $C^1$ function defined on a neighborhood of $b$ such that for $t\rightarrow b^- $, $g_2(t) >tg_2'(t)$ and 
\begin{equation}\label{condipourg2}
h(t)\leq \lambda_\mu \frac{t^2(1-g_2(t))}{g_2(t)-tg_2'(t)},\quad t\rightarrow b^-
\end{equation}
If $\underset{t\rightarrow b}{\lim}\,v(t)=0 $, then $$1-q(t)\geq \frac{C_P(\mu)}{Z}v(t)\frac{g_2(t)}{t}, \quad t\rightarrow b^- . $$
\end{lemm}

Under a growth control on $h$, the condition $v(t) \underset{x\rightarrow a}{\rightarrow} 0$ can easily be verified. One needs however to distinguish the case where $a$ is finite from the one where it is infinite. For $a=-\infty$: if $h(t)\leq c|t|^\alpha$ for $t\rightarrow -\infty$ with $c>0$ and $\alpha\leq 2$, then for $x\rightarrow -\infty$, $$-\lambda_\mu\int_0^x\frac{u}{h(u)}\,du \leq -\frac{\lambda_\mu }{c}\int_x^0\frac{1}{|u|^{\alpha-1}} du \underset{x\rightarrow-\infty}{\rightarrow} -\infty,$$ and hence $v(t) \underset{x\rightarrow -\infty}{\rightarrow} 0$. For $a>-\infty$: if $h(t)\leq c(t-a)^\alpha$ for $t\sim a$ with $c>0$ and $\alpha\geq 1$, then for $x\in ]a,0[$: $$-\lambda_\mu\int_{0}^x\frac{u}{h(u)}\,du \leq \lambda_\mu \int_x^{0}\frac{u}{(u-a)^\alpha}du = \lambda_\mu \int_x^0\frac{du}{(u-a)^{\alpha-1}} + a\lambda_\mu \int_x^0 \frac{du}{(u-a)^\alpha} \underset{x\rightarrow a}{\rightarrow} -\infty,$$ and hence $v(t) \underset{x\rightarrow a}{\rightarrow} 0$. The same results are valid at $b$, only replacing $(t-a)^\alpha$ by $(b-t)^\alpha$ when $b<+\infty$.

\begin{prop}\label{uncafinipasfini}
If $a=-\infty$ and $c_1 |t|^\beta\leq h(t) \leq c_2 |t|^{\alpha}$ for $t\rightarrow -\infty$ with $c_1,\,c_2>0$, $\alpha\leq2$, and $\beta\in [2\alpha-2,\,\alpha] $, then $q$ satisfies ($\ref{g_1}$) with $g_1(t)=1-\frac{c_2}{\lambda_\mu} |t|^{\alpha-2}$. Moreover, ($\ref{bornitude} $) is also satisfied. \\
If $a>-\infty$ and $c_1(t-a)^\beta h(t)\leq c_2(t-a)^\alpha$ for $t\rightarrow a^+$ with $c_1,\,c_2>0$, $\alpha\geq 1$ and $\beta\leq 2\alpha $, then $q$ satisfies ($\ref{g_1}$) with $g_1(t)=1-\frac{4c_2}{\lambda_\mu a^2} (t-a)^{\alpha}$.
\end{prop}

\begin{dem}
Let us compute, if $a=-\infty$ and $t<0$: $$\lambda_\mu\frac{t^2(1-g_1(t))}{g_1(t)-tg_1'(t)} = \frac{c_2t^2|t|^{\alpha-2}}{1-\frac{c_2}{\lambda_\mu}(3-\alpha)|t|^{\alpha-2}} \geq c _2|t|^{\alpha},$$ for $t\sim -\infty$. 
If $a>-\infty$ and $t\in]a,\frac{a}{2}[$, $$\lambda_\mu\frac{t^2(1-g_1(t))}{g_1(t)-tg_1'(t)} = \frac{4c_2}{a^2}\frac{t^2(t-a)^{\alpha}}{1+\frac{4c_2}{\lambda_\mu a^2}(t-a)^{\alpha-1}[\alpha t-(t-a)]} \geq c_2(t-a)^\alpha ,$$ since $\alpha t-(t-a)<0$ and $t^2\geq \frac{a^2}{4}$.

Hence Lemma $\ref{lemmcondipourg1}$ and the above remark give us that $q$ satisfies $(\ref{g_1})$ in both case.\\
In regards to $(\ref{bornitude})$, in the first case, one has $$\left|\frac{x(1-g_1(x))}{\sqrt{h(x)}}\right|=\frac{c_2}{\lambda_\mu\sqrt{h(x)}}|x|^{\alpha-1}\leq \frac{c_2}{\lambda_\mu\sqrt{c_1}}|x|^{\alpha-1-\frac{\beta}{2}} $$ which is bounded at $-\infty$ because $\alpha-1-\frac{\beta}{2}\leq \alpha-1-\frac{2\alpha-2}{2}=0$. Similarly in the second case,
$$ \left|\frac{x(1-g_1(x))}{\sqrt{h(x)}}\right|= \frac{4c_2\,|x|}{\lambda_\mu \,a^2}\cdot\frac{|x-a|^\alpha}{\sqrt{h(x)}}\leq \frac{4c_2}{\lambda_\mu\, |a|}\frac{1}{\sqrt{c_1}}|x-a|^{\alpha-\frac{\beta}{2}} $$ which is bounded at $a>-\infty$ because $\beta\leq 2\alpha$.
\end{dem}
The same result is valid at $b$, using Lemma $\ref{lemmcondipourg2}$. We summarize it in the following Proposition:

\begin{prop}\label{resumechfini}
Assume that one of these two conditions is verified at $a$: 
\begin{itemize}
\item either $a=-\infty$ and $c_1 |t|^{2\alpha-2}\leq h(t)\leq c_2 |t|^\alpha$ for $t\rightarrow -\infty$ with $\alpha\leq 2 $ and $c_1,\,c_2>0$,
\item or $a>-\infty$ and $c_1 (t-a)^2\leq h(t)\leq c_2(t-a)$ for $t\rightarrow a^+$ with $c_1,\,c_2>0$,
\end{itemize}
and one of these two conditions is satisfied at $b$:
\begin{itemize}
\item either $b=+\infty$ and $c_1 t^{2\alpha-2}\leq h(t)\leq c_2 t^\alpha$ for $t\rightarrow +\infty$ with $\alpha\leq 2 $ and $c_1,\,c_2>0$,
\item or $b<+\infty$ and $c_1 (b-t)^2\leq h(t)\leq c_2(b-t)$ for $t\rightarrow b^-$ with $c_1,\,c_2>0$.
\end{itemize}
Then the constant $C_h$ defined in Theorem $\ref{thmcentr}$ is finite.
\end{prop}

\subsection{Application: stability of the $\Gamma(s,\theta)$ distributions, $s>0$, $\theta>0$}\label{gammadistrib}

Let us consider the $\Gamma(s,\theta)$ distribution $d\mu(x) := \frac{x^{s-1}e^{-\frac{x}{\theta}}}{\Gamma(s)\theta^s}\textbf{1}_{\mathbb{R}_+}(x) dx$, where $\Gamma$ denotes the Euler's Gamma function. This probability measure is the reversible measure of the Laguerre process with following generator on $[0,\infty)$:
$$Lf(x) = xf''(x) + (s-\frac{1}{\theta} x)f'(x) .$$
We have $\Gamma(f)=x(f')^2$, $C_P(\mu)=\theta$ and the first Laguerre polynomial $f_0(x)=\frac{s-\frac{x}{\theta}}{\sqrt{s}}$ is a normalized eigenfunction attaining the spectral gap. Hence $I=]-\infty,\sqrt{s}]$, so we get $h(x)=\theta^{-1}\left( s-\sqrt{s}\,x\right)$. So Assuptions $\ref{spectralgapattained}$, $\ref{gammaf0measuable}$ and $\ref{criticalmaxima}$ are satisfied. Moreover, one can see that we are also under the conditions of Proposition $\ref{resumechfini}$, with at the infimum, $a=-\infty$ and $c_1 \leq h(t)\leq c_2 |t|$ for  $c_1=\frac{s}{\theta} $, $c_2=2\frac{\sqrt{s}}{\theta}$ (and $\alpha=1$), and at the supremum, $b=\sqrt{s}<+\infty$, $c_1 (b-t)\leq h(t)\leq c_2(b-t)$ for $c_1=c_2=\frac{\sqrt{s}}{\theta}$.
In this case the three conditions ($\ref{hypnu}$) are reduced to only two:
\begin{equation}\label{normalizationreduced}
\int x\,d\nu = s\theta \quad \mathrm{and} \quad \int x^2\, d\nu = s(s+1)\theta^2
\end{equation}
 meaning that we only ask $\nu$ to have the same first and second moments as $\mu$.
So we can apply Theorem $\ref{thmcentralexploitable}$ for the gamma distributions $\Gamma(s,\theta)$, with $s,\theta>0$: for all measures $\nu$ on $\mathbb{R}_+$ normalized as in ($\ref{normalizationreduced}$) and satisfying the following Poincaré inequality with sharp constant $C_P(\nu)$: 
\begin{equation}\label{poincargamm}
\mathrm{Var}_\nu (f) \leq C_P(\nu)\int_0^{+\infty} xf'(x)^2d\nu(x), 
\end{equation} it holds for a finite constant $C_h$:
\begin{equation}\label{stabgamm}
W_1(\mu^*,\nu^*) \leq C_h\left(\sqrt{\delta} + \sqrt{C_P(\nu)}\,\delta \right), 
\end{equation}

where $W_1$ is the $1$-Wasserstein distance (see Definition $\ref{w1int}$), $\delta:= \frac{ C_P(\nu) - \theta}{C_P(\nu)}$, $\nu^*:=f^{\#}_0 \nu$, and $$d\mu^*(x)=\left(\frac{\sqrt{s}}{e}\right)^s\frac{1}{\Gamma(s)}(\sqrt{s}-x)^{s-1}e^{\sqrt{s}\,x}\mathbf{1}_{]-\infty,\sqrt{s}]}(x)dx. $$ 

\subsection{Application : stability for Brascamp-Lieb inequalities}\label{brascampliebstab}

Let $M=\mathbb{R}$ and $d\mu=e^{-\phi}dx $ be a probability measure with $\phi:\mathbb{R}\rightarrow \mathbb{R}$ strictly convex and $C^2$. Then $\mu$ satisfies the following inequality introduced by Brascamp and Lieb in 1976 in \cite{brascamplieb} : $$ \mathrm{Var}_{\mu}(f) \leq \int \frac{f'^2}{\phi''}d\mu $$ for all functions $f$ in the weighted Sobolev space $H^{1,2}_{\phi,\mu}= \{g \in L^2(\mu)| \frac{g'}{\sqrt{\phi''}}\in L^2(\mu)\}$.
Then $$Lf = \frac{1}{\phi''}f'' - \left(\frac{\phi'}{\phi''} + \frac{\phi^{(3)}}{\phi''^2}\right)f'$$ is a Markov diffusion generator with reversible measure $\mu$ satisfying a Poincaré inequality with sharp constant $1$ and carré du champ operator $\Gamma(f)=\frac{f'^2}{\phi''}$. Moreover $$L(\phi')=-\phi', $$ hence the spectral gap is attained with $f_0=\phi'$ which is centered. One easily sees that $\phi'$ is a bijection, so $I=\mathbb{R}$ and Assumption $\ref{gammaf0measuable}$ is satisfied with $h=\phi''\circ\phi'^{-1}$. 
The push forward measure $\mu^*$ is then the moment measure of $\phi$ (see \cite{momentmeas,santa}) and has density $$d\mu^* = \frac{\exp\left( -\phi(\phi'^{-1}(t)) \right)}{\phi''(\phi'^{-1}(t))}\, dt. $$ 
In this setting, we can get two stability results with respect to two different distances. Firstly, Theorem $\ref{stabilityhminore}$ says that if $\phi'' \geq \kappa >0$ (meaning that $\phi$ is uniformly convex) then $\mu^*$ satisfies a stability result in total variation distance. It is a direct generalization of the Gaussian case \cite{U}.

\begin{corol}
Let $M=\mathbb{R}$ and $d\mu=e^{-\phi}\,dx $ be a $C^2$ uniformly log-concave probability measure, i.e. $\forall x\in\mathbb{R},\, \phi''(x)\geq \kappa>0$. Let $\nu$ be a probability measure on $\mathbb{R}$ verifying the following Poincaré inequality: for all functions $f$ in the weighted Sobolev space $H^{1,2}_{\phi,\nu}= \{g \in L^2(\nu)| \frac{g'}{\sqrt{\phi''}}\in L^2(\nu)\}$, $$ \mathrm{Var}_{\nu}(f) \leq C_P(\nu)\int \frac{f'^2}{\phi''}d\nu, $$ and such that $\int \phi' \,d\nu = 0 $, $\int \phi'^2\,d\nu=1 $ and $\int \phi''\,d\nu \leq \int \phi''\,d\mu$. So $C_P(\nu)\geq 1$, and moreover $\delta:= \frac{ C_P(\nu) - 1}{C_P(\nu)}$ satisfies $$d_{TV}(\mu^*,\nu^*) \leq \frac{4}{\sqrt{\kappa}}\left(\sqrt{\delta} + \sqrt{C_P(\nu)}\,\delta \right), $$ where $\mu^*=\phi'_{\#}\mu$ and $\nu^*=\phi'_{\#}\nu $.
\end{corol}

Secondly, Theorem $\ref{thmcentralexploitable}$ yields a stability in $1$-Wasserstein distance under a growth control of the second derivative of the potential $\phi$ at infinity.
\begin{corol}
Let $M=\mathbb{R}$ and $d\mu=e^{-\phi}\,dx $ be a $C^2$ strictly log-concave probability measure. Let $\nu$ be a probability measure on $\mathbb{R}$ verifying the following Poincaré inequality: for all functions $f$ in the weighted Sobolev space $H^{1,2}_{\phi}= \{g \in L^2(\nu)| \frac{g'}{\sqrt{\phi''}}\in L^2(\nu)\}$, $$ \mathrm{Var}_{\nu}(f) \leq C_P(\nu)\int \frac{f'^2}{\phi''}d\nu, $$ and such that $\int \phi' \,d\nu = 0 $, $\int \phi'^2\,d\nu=1 $ and $\int \phi''\,d\nu \leq \int \phi''\,d\mu$. So $C_P(\nu)\geq 1$, and if for some $\alpha\leq 2$, $$c_1 |t|^{2\alpha-2}\leq \phi''(t) \leq c_2 |t|^\alpha, \quad t\rightarrow \pm \infty,$$ then
$\delta:= \frac{ C_P(\nu) - 1}{C_P(\nu)}$ satisfies $$W_1(\mu^*,\nu^*) \leq C_h\left(\sqrt{\delta} + \sqrt{C_P(\nu)}\,\delta \right), $$
where $\mu^*=\phi'_{\#}\mu$ and $\nu^*=\phi'_{\#}\nu $ and $C_h$ is a finite constant defined by ($\ref{defdeCh}$).
\end{corol}

\subsection{Application : stability of the Poincaré constant on the sphere $\mathbb{S}^d$}

Let $M=\mathbb{S}^d$, $d\geq 1$ and $\mu=\frac{\Gamma(\frac{d}{2})}{2\pi^{\frac{d}{2}}} \mathrm{vol}_{\mathrm{S}^d}$ the riemannian volume measure on the sphere, normalized to be a probability measure, where $\Gamma$ denotes Euler's Gamma function. Then $\mu$ is the reversible probability measure of the spherical Laplacian generating the Brownian motion on the sphere, which can be written as the restriction to the sphere of $$Lf(x)= \sum_{i,j=1}^{d+1}\left(\delta_{ij}-x_ix_j\right)\partial^2_{ij}f - d\sum_{i=1}^{d+1}x_i\partial_if, $$ if $f$ is the restriction to the sphere of a smooth function on a neighborhood of the sphere in $\mathbb{R}^{d+1}$ (see \cite[Section 2.2.1]{BGL} for more details), and where $\delta_{ij}$ denotes the Kronecker delta. 

We have $\Gamma(f)=  \sum_{i,j=1}^{d+1}\left(\delta_{ij}-x_ix_j\right)\partial_{i}f\,\partial_jf$, $C_P(\mu)=\frac{1}{d}$ and all coordinate functions $x\mapsto x_j$ restricted to the sphere are centered eigenfunctions attaining the spectral gap. Let us choose $f_0(x)=\sqrt{d+1}\, x_1$ which has variance one. Hence $I=[-\sqrt{d+1},\sqrt{d+1}]$ and $\Gamma(f_0)=(d+1)\left(1-x_1^2\right)$, so we get $h(t)=(d+1)-t^2 $. So Assuptions $\ref{spectralgapattained}$, $\ref{gammaf0measuable}$ and $\ref{criticalmaxima}$ are satisfied. Furthermore, at the infimum $a=-\sqrt{d+1}$, one sees $\sqrt{d+1}\left(\sqrt{d+1}+t\right)\leq h(t)\leq 2\sqrt{d+1} \left(\sqrt{d+1}+t\right) $, and at the supremum $b=\sqrt{d+1}$, one similarly sees $\sqrt{d+1}\left(\sqrt{d+1}-t\right)\leq h(t)\leq 2\sqrt{d+1} \left(\sqrt{d+1}-t\right) $. Hence the conditions of Theorem $\ref{thmcentralexploitable}$ are satisfied. In this case the three conditions ($\ref{hypnu}$) are reduced to only two:
\begin{equation}\label{normalizationreducedforsphere}
\int_{\mathbb{S}^d} x_1\,d\nu = 0 \quad \mathrm{and} \quad \int_{\mathbb{S}^d} x_1^2\, d\nu = \frac{1}{d+1}.
\end{equation}
So Theorem $\ref{thmcentralexploitable}$ applies to the canonical probability distribution on $\mathbb{S}^d$.

\begin{corol}
Let $M=\mathbb{S}^d$, $d\geq 1$ and $\mu=\frac{\Gamma(\frac{d}{2})}{2\pi^{\frac{d}{2}}} \mathrm{vol}_{\mathrm{S}^d}$. Let $\nu$ be a probability measure on $\mathbb{S}^d$ normalized as in ($\ref{normalizationreducedforsphere}$) and satisfying the following Poincaré inequality with sharp constant $C_P(\nu)$: for all functions $f$ which is the restriction to the sphere of a $C^2$ function on $\mathbb{R}^{d+1}$,
\begin{equation}\label{poincarspher}
\mathrm{Var}_\nu (f) \leq C_P(\nu)\int_{\mathbb{S}^d} \sum_{i,j=1}^{d+1}\left(\delta_{ij}-x_ix_j\right)\partial_{i}f\,\partial_jf\,d\nu(x).
\end{equation} Then it holds for a finite constant $C_h$,
\begin{equation}\label{stabspher}
W_1(\mu^*,\nu^*) \leq C_h\left(\sqrt{\delta} + \sqrt{C_P(\nu)}\,\delta \right), 
\end{equation}
where $W_1$ is the $1$-Wasserstein distance (see Definition $\ref{w1int}$), $\delta:= \frac{ C_P(\nu) - d}{C_P(\nu)}$, $\nu^*:=x_1^{\#} \nu$, and $$d\mu^*(t)=\frac{1}{Z}\left(d+1-t^2 \right)^{\frac{d}{2}-1}\mathbf{1}_{[-\sqrt{d+1},\sqrt{d+1}]} dx, $$  with $Z$ a normalization constant.
\end{corol}

\section{Comparison}\label{comparison}

A natural question is to compare the scope of application of Theorem $\ref{stabilityhminore}$ and Theorem $\ref{thmcentralexploitable}$. Indeed, Theorem $\ref{stabilityhminore}$ require $h$ to be uniformly bounded by below, while Theorem $\ref{thmcentralexploitable}$ only require a growth control on $h$ at the boundary of $I$. In this section we get stability results involving the Kolmogorov distance under the assumptions of Theorems $\ref{stabilityhminore}$ and $\ref{thmcentralexploitable}$. Afterwards, we will see that the order of the bound obtained from Theorem $\ref{stabilityhminore}$ is stronger. Let us begin recalling:

\begin{prop}\label{kolmoinfw1} \cite[Prop 1.2]{fundamStein}
If $\mu^*$ has bounded density $\rho$ with respect to the Lebesgue measure on $\mathbb{R}$, then $$d_K(\nu^*,\mu^*)\leq 2\sqrt{C\, W_1(\nu^*,\mu^*)} $$ where $\rho\leq C$, $W_1$ is the $1$-Wasserstein distance, and $d_K(\nu^*,\mu^*)=\underset{t\in\mathbb{R}}{\sup}\left|\int_{-\infty}^t d\nu^* - \int_{-\infty}^t d\mu^* \right|$ is the Kolmogorov distance.
\end{prop}

Combining Proposition $\ref{kolmoinfw1}$ above with Theorem $\ref{thmcentralexploitable}$ gives the following bound on the Kolmogorov distance.

\begin{corol}\label{kolmoavecCh}
Let $(L,\mu)$ and $h$ statisfy the requirements of Theorem $\ref{thmcentralexploitable}$. If moreover $\frac{v(x)}{h(x)}$ is bounded by a constant $C$, then it holds for $\delta:= \frac{ C_P(\nu) - C_P(\mu)}{C_P(\nu)\,C_p(\mu)}$, $$d_{K}(\mu^*,\nu^*)\leq 2\sqrt{C\,C_h}\sqrt{\sqrt{\delta} + \sqrt{C_P(\nu)}\,\delta }, $$ for all measures $\nu$ normalized as in ($\ref{hypnu}$) and satisfying a Poincaré inequality with sharp constant $C_P(\nu)$.
\end{corol}

On the other hand, it is well-known that the total variation distance controls the Kolmogorov one: $d_K(\nu^*,\mu^*)\leq d_{TV}(\nu^*,\mu^*) $. Hence, Theorem $\ref{stabilityhminore}$ implies the following.

\begin{corol}\label{kolmoavechminoree}
Let $(L,\mu)$ and $h\geq \kappa>0$ statisfy the requirements of Theorem $\ref{stabilityhminore}$. Then it holds for $\delta:= \frac{ C_P(\nu) - C_P(\mu)}{C_P(\nu)\,C_p(\mu)}$, $$d_{K}(\mu^*,\nu^*)\leq \frac{4}{\sqrt{\kappa}}\left(\sqrt{\delta} + \sqrt{C_P(\nu)}\,\delta \right), $$ for all measures $\nu$ normalized as in ($\ref{hypnu}$) and satisfying a Poincaré inequality with sharp constant $C_P(\nu)$.
\end{corol}

In case where $(L,\mu)$ satisfies both assumptions of Corollary $\ref{kolmoavecCh}$ and Corollary $\ref{kolmoavechminoree}$, one can see that the bound of Corollary $\ref{kolmoavechminoree}$ is better because it is of order $\frac{1}{2}$ against $\frac{1}{4}$ for the bound of Corollary $\ref{kolmoavecCh}$. However, in the example of the gamma distributions, while Corollary $\ref{kolmoavechminoree}$ does not apply, although non-optimal, Corollary $\ref{kolmoavecCh}$ gives us the following stability result in Kolmogorov distance.

\begin{corol}
For all $\theta,s>0$ and for all measures $\nu$ on $\mathbb{R}_+$ normalized as in ($\ref{normalizationreduced}$) and satisfying the Poincaré inequality ($\ref{poincargamm}$) with sharp constant $C_P(\nu)$, it holds: $$d_{K}(\mu^*,\nu^*)\leq 2\sqrt{C\,C_h}\sqrt{\sqrt{\delta} + \sqrt{C_P(\nu)}\,\delta }, $$ 
where $\delta:= \frac{ C_P(\nu) - \theta}{C_P(\nu)}$, $C=\left( \sqrt{s}\right)^{s+\sqrt{s}}\frac{e^{1-s}}{\Gamma(s)} $, $\nu^*:=\left( \sqrt{s}-\frac{x}{\theta\sqrt{s}}\right)^{\#} \nu$, and $$d\mu^*(x)=\left(\frac{\sqrt{s}}{e}\right)^s\frac{1}{\Gamma(s)}(\sqrt{s}-x)^{s-1}e^{\sqrt{s}\,x}\mathbf{1}_{]-\infty,\sqrt{s}]}(x)dx. $$ 
\end{corol}

\begin{acknowledgments}
I would like to thank my PhD advisors Franck Barthe and Max Fathi for their availability and all their fruitful advice.
\end{acknowledgments}

\bibliographystyle{plain}
\bibliography{mabibliographie}

\textbf{Contact information:} Institut de Mathématiques de Toulouse\\
E-mail: jordan.serres@math.univ-toulouse.fr

\end{document}